\documentclass{amsart}
\pdfoutput=1
\usepackage{amsmath, amsthm, amssymb}
\usepackage{hyperref} 

\DeclareMathOperator{\Frob}{Frob} 
\DeclareMathOperator{\Gal}{Gal} 
\DeclareMathOperator{\Mass}{Mass} 
\DeclareMathOperator{\Res}{Res} 
\DeclareMathOperator{\Sym}{Sym} 
\DeclareMathOperator{\Norm}{Norm} 
\DeclareMathOperator{\GL}{GL}

\newtheorem{theorem}{Theorem}[section]
\newtheorem{lemma}[theorem]{Lemma}
\newtheorem{definition}[theorem]{Definition}
\newtheorem{proposition}[theorem]{Proposition}
\newtheorem*{algorithmA}{Algorithm A}

\theoremstyle{remark}
\newtheorem{remark}[theorem]{Remark}

\begin{document}

\title
	[Definite unitary automorphic forms]
	{Explicit calculations of automorphic forms for definite unitary groups}

\author{David Loeffler}
\email{david.loeffler05@imperial.ac.uk}
\address{Department of Mathematics \\ Imperial College \\ South Kensington \\ London SW7 2AZ}


\date{22 January 2008, revised 11 August 2008}

\begin{abstract}
I give an algorithm for computing the full space of automorphic forms for definite unitary groups over $\mathbb{Q}$, and apply this to calculate the automorphic forms of level $G(\hat{\mathbb{Z}})$ and various small weights for an example of a rank 3 unitary group. This leads to some examples of various types of endoscopic lifting from automorphic forms for $U_1 \times U_1 \times U_1$ and $U_1 \times U_2$, and to an example of a non-endoscopic form of weight $(3,3)$ corresponding to a family of 3-dimensional irreducible $\ell$-adic Galois representations. I also compute the $2$-adic slopes of some automorphic forms with level structure at 2, giving evidence for the local constancy of the slopes.
\end{abstract}

\maketitle
\tableofcontents

\section{Introduction}

If $G$ is a reductive algebraic group over $\mathbb{Q}$, then the spaces of automorphic forms for $G$ (of a given level and infinity-type) are well-known to be finite-dimensional vector spaces with an action of the Hecke algebra. However, for the vast majority of groups it is not known how to actually calculate the dimensions of these spaces and how the Hecke operators act. For $\GL_2$, the case of classical modular forms, there are well-known algorithms based on modular symbols (see e.g. \cite{stein}), but in general little is known.

One case in which calculations are possible is where the group $G$ satisfies the condition that all arithmetic subgroups are finite. In this case, Gross has shown \cite{gross-algebraic} that the theory of automorphic forms can be set up entirely algebraically, without analytic hypotheses. As has been observed by various authors \cite{lansky-pollack, dembele-hilbert}, Gross's spaces are (at least in theory) computable. In this article, I shall apply this to an example of a definite unitary group over $\mathbb{Q}$, and show how the full space of automorphic forms may be calculated for various small weights and levels. We do not need to assume that the class number of $G$ is 1; we do assume that the underlying quadratic extension of $\mathbb{Q}$ has class number 1, but this is more a calculational convenience than a fudamental limitation of the algorithm.

The results of these computations give rise to some interesting specimens of endoscopic liftings of automorphic forms arising from $U_2 \times U_1$ and $U_1 \times U_1 \times U_1$, and also examples of non-endoscopic forms, which give nontrivial families of 3-dimensional $\ell$-adic Galois representations.

The structure of this paper is as follows. Section \ref{sect:algorithm} explains the algorithms used to determine the class number of $G$ and the coset decomposition of the Hecke operators; these depend on the choice of level group $K$ but not the weight. In section \ref{sect:weights}, we use the results of these computations to calculate the spaces of forms of various small weights. In section \ref{sect:higher}, we shall see how to generalise this to non-maximal level groups, and present some examples illustrating the 2-adic continuity of the Hecke eigenvalues. In the final section, we discuss some implementation details.

\section{Outline of the algorithm}
\label{sect:algorithm}

\subsection{Definitions}

Recall that if $E$ is an imaginary quadratic number field and $n \ge 1$, there is a corresponding definite unitary group $G = U_{n,E}$, which is the group scheme over $\mathbb{Z}$ defined by 
\[ G(A) = \{g \in \GL_n(A \otimes_\mathbb{Z} \mathcal{O}_E) \ |\ g g^\dagger = 1\},\]
where $A$ is an arbitrary commutative ring (and $g^\dagger$ is the conjugate transpose of $g$).

The following properties of $G$ are elementary:

\begin{proposition}
$G(\mathbb{R})$ is the usual real unitary group $U_n$, which is connected and compact. If $p$ is any prime which splits completely in $E$, then $G$ is split at $p$, and moreover $G(\mathbb{Q}_p)$ is isomorphic to $\GL_n(\mathbb{Q}_p)$. There are two such isomorphisms, corresponding to the two embeddings $E \hookrightarrow \mathbb{Q}_p$, and they are interchanged by the inverse transpose involution of $\GL_n$. 
The centre of $G$ is the norm torus $U(1) = \Res_{E/\mathbb{Q}}^{(1)}(\mathbb{G}_m)$.
\end{proposition}

Since $G$ is compact at infinity, it certainly satisfies the conditions of \cite{gross-algebraic}; so the space of automorphic forms for $G$ of level $K$ (where $K$ is an open compact subgroup of $G(\mathbb{A}_f)$) and weight $V$ (where $V$ is an irreducible algebraic representation of $G$ over any field $F$ of characteristic 0) is the $F$-vector space of functions $f: G(\mathbb{A}_f) \to V$ such that 
\[f(\gamma g k) = \gamma \circ f(g)\quad \forall \gamma \in G(\mathbb{Q}),\ k \in K. \]

For the rest of this paper, we shall fix $n = 3$; however, the methods can clearly be applied to arbitrary $n$, and indeed to general totally definite Hermitian spaces over $\mathcal{O}_E$ (of any integral equivalence class). We shall, however, make the simplifying assumption that $E$ has class number 1. Initially, we shall take $K = G(\hat{\mathbb{Z}}) = \prod_\ell G(\mathbb{Z}_\ell)$, so we obtain the automorphic forms of ``level 1'' in some sense. 

\begin{remark}
If $p$ is not ramified in $E$, the subgroup $K$ is a hyperspecial maximal compact subgroup of $G(\mathbb{Q}_p)$. If $p$ is ramified, then $K$ is special maximal but not hyperspecial, unless $p = 2$ in which case it is not even maximal \cite[\S1.10]{rogawski}.
\end{remark}

\begin{theorem}
The irreducible algebraic representations of $G$ over $E$ (or any field containing it) are parametrised by triples of integers $(a, b, c)$ with $a, b \ge 0$, where the representation $V_{a,b,c}$ is the unique highest weight direct summand of $W_{a,b,c} = \Sym^a V \otimes \Sym^b V^* \otimes \det^c$  (where $V = V_{1,0,0}$ is the standard representation). The central character of $V_{a,b,c}$ is $\begin{pmatrix}z & & \\ & z & \\ & & z\end{pmatrix} \mapsto z^{a - b + 3c}$.
\end{theorem}

\begin{proof}
Since $G$ is isomorphic to $\GL_3$ over $E$, this follows from standard results on the representation theory of $\GL_3$, see e.g.~\cite[\S 15.5]{fulton-harris}.
\end{proof}

\begin{remark}
We can identify $V_{a,b,c}$ explicitly: there is a natural contraction map $W_{a,b,c} \to W_{a-1,b-1,c}$, and $V_{a,b,c}$ is the kernel of this map. We shall not use this remark, however, as in the calculations below we will always have already calculated automorphic forms of weight $(a-1,b-1,c)$ when we come to deal with those of weight $(a,b,c)$; thus we can easily identify forms which we have seen before.
\end{remark}

\subsection{\texorpdfstring{Hecke operators and $K$-classes}{Hecke operators and K-classes}}

It is well known that the set of double cosets
\[ G(\mathbb{Q}) \backslash G(\mathbb{A}_f) / K\]
is finite. We shall refer to these as {\bf $K$-classes}. As shown in \cite{gross-algebraic}, if $\mu_1, \dots, \mu_r$ is a set of representatives for the $K$-classes, then the space of automorphic forms for $G$ of level $K$ and weight $V$ is isomorphic to 
\[ \bigoplus_{i=1}^r V^{\Gamma_i},\]
where $\Gamma_i = G(\mathbb{Q}) \cap \mu_i K \mu_i^{-1}$, via the map $f \mapsto \left(f(\mu_1), \dots, f(\mu_r)\right)$. Note that the groups $\Gamma_g$ are arithmetic subgroups of $G$, and are thus finite.

If we know a set of class representatives $\mu_i$ and the associated groups $\Gamma_i$, then it is elementary to calculate the $\Gamma_i$-invariants in $V_{a,b,c}$ for each $i$ and thus read off the dimension of the space of automorphic forms. 

The other piece of information we are interested in is the action of the Hecke algebra. This is the commutative algebra generated by double cosets $K g K$, with the action of such a coset on automorphic forms being given by 
\[ \left([K g K] f\right)(x) = \sum_j f(x g_j)\]
where $K g K = \bigsqcup g_j K$. 

So we must decompose $K g K$ into single cosets $g_j K$. For each of these, we could choose a coset representative in the form $\gamma_j \mu_{c(j)}$, for some $c(j) \in \{1 \dots r\}$ and $\gamma_j \in G(\mathbb{Q})$. This would immediately allow us to calculate $\left([K g K] f\right)(1)$, if $f$ is an automorphic form given as an $r$-tuple $f(\mu_1), \dots f(\mu_r)$, since 
\[\left([K g K] f\right)(1) = \sum_j f(\gamma_j \mu_{c(j)}) = \sum_j \gamma_j \circ f(\mu_{c(j)}).\]

What we actually want, however, is $\left([K g K] f\right)(\mu_i)$; so for each $i \in \{1, \dots, r\}$ we need to find elements $\gamma_{ij} \in G(\mathbb{Q})$ and $c(i,j) \in \{1, \dots, r\}$ such that 
\[\mu_i K g K = \bigsqcup_j \gamma_{ij} \mu_{c(i,j)} K.\]

For any given element $g \in G$, finding these quantities is a finite search. In fact, one has:

\begin{algorithmA}
Let $r$, $s$, $t$ be given elements of $G(\mathbb{A}_f)$, which are integral outside a finite set $S$ of primes split in $E$. Then the following algorithm finds all elements of $G(\mathbb{Q}) \cap r K s K t$:
\begin{enumerate}
\item Find some $\lambda \in \mathcal{O}_E$ such that $\lambda r s t$, regarded as an element of $M_3(\mathbb{A}_{f,E})$, is integral at all places of $E$.
\item Enumerate all matrices $g \in M_3(\mathcal{O}_E)$ such that $g g^\dag = \Norm(\lambda)$.
\item For each $g$ in the above list, set $\gamma = \lambda^{-1} g$ and calculate the elementary factors of $r^{-1} \gamma t^{-1}$ at each prime in $S$. If these coincide with those of $s$, output $\gamma$.
\end{enumerate}
\end{algorithmA}

Calculating the elementary factors of a $p$-adic matrix is easy, but in practice one can often short-cut this step, since if $\nu_p \lambda$ is small then there are very few possibilities for the elementary factors and these can be distinguished by calculating the determinant. The most computationally intensive step is usually (2). For convenience, let us define an {\it $r$-good matrix} to be an element of $M_3(\mathcal{O}_E)$ such that $g g^\dag = r$. Enumerating all the $r$-good matrices for a given $r$ is clearly a finite check, but possibly a long one if $r$ is large -- an algorithm using lattices is described in \S \ref{sect:rgood} below.

As is well known, the Hecke algebra factorises into a restricted tensor product of local Hecke algebras at each prime; if $p$ is a prime split in $E$, and $\mathfrak{p}$ is a choice of prime above $p$ in $E$, then the local Hecke algebra of $K$ at $p$ is generated by the double cosets $K \eta K$, where $\eta$ is one of the elements
\begin{align*}
 \eta_{\mathfrak{p}, 1} &= \begin{pmatrix} 1 &  & \\ & 1 & \\ & & \mathfrak{p}/\overline{\mathfrak{p}} \end{pmatrix}, & 
 \eta_{\mathfrak{p}, 2} &= \begin{pmatrix} 1 &  & \\ & \mathfrak{p}/\overline{\mathfrak{p}} & \\ & & \mathfrak{p}/\overline{\mathfrak{p}} \end{pmatrix}, &
 \eta_{\mathfrak{p}, 3} &= \begin{pmatrix} \mathfrak{p}/\overline{\mathfrak{p}} &  & \\ & \mathfrak{p}/\overline{\mathfrak{p}} & \\ & & \mathfrak{p}/\overline{\mathfrak{p}} \end{pmatrix}.
\end{align*}

These are the operators to which we shall apply the above method. We shall write $T_{\mathfrak{p},i}$ for the operator on automorphic forms corresponding to the double coset $K \eta_{\mathfrak{p},i} K$. Note that $T_{\mathfrak{p}, 3}$ is central, and acts on weight $(a,b,c)$ forms as scalar multiplication by $(\mathfrak{p}/\overline{\mathfrak{p}})^{a-b+3c}$.

\subsection{A refinement}

The above discussion assumes that the class representatives $\mu_i$ are known. However, we can use a sort of bootstrapping approach to find the class representatives at the same time as the Hecke operators.

\begin{lemma}\label{lemma:localbound}
The double cosets $K \eta_{\mathfrak{p}, 1} K$ and $K \eta_{\mathfrak{p},2} K$ both contain $p^2 + p + 1$ single cosets.
\end{lemma}

\begin{proof} This is a purely local computation. We shall consider $\eta = \eta_{\mathfrak{p}, 1}$; the other case is very similar. Let $H = G(\mathbb{Z}_p)$; then $H$ acts transitively by right multiplication on nonzero row vectors over $\mathbb{F}_p$, and $H \cap \eta H \eta^{-1}$ is the stabilizer of the line $(0,0,*)$. So this subgroup has index $|\mathbb{P}^2(\mathbb{F}_p)| = p^2 + p + 1$; and thus $H \eta H$ should decompose into $p^2 + p + 1$ single cosets (which can easily be written down explicitly).
\end{proof}

Now, let us assume that we know a set of representatives for a subset of the class set; we can always start with just the principal class $G(\mathbb{Q}) K$, represented by the element 1. Let us choose a prime $p$. Using Algorithm A above, we can find all of the single cosets contained in $K \eta_{\mathfrak{p}, 1} K$ which have a representative of the form $\gamma \mu$ with $\gamma \in G(\mathbb{Q})$ and $\mu$ one of the subset of class representatives we know about.

If we have found less than $p^2 + p + 1$ cosets, then the lemma shows that we have not found the full class set, and moreover, taking a local coset representative gives an explicit element $\mu$ of $G(\mathbb{A}_f)$ (supported at $p$) which is not in any of the $K$-classes previously found; so we can add this to the list, calculate the associated group $\Gamma$ (by applying Algorithm A again, with $s = 1$) and repeat.

Clearly, we need a criterion which will allow us to determine when we have found the complete class set. To do this, we shall introduce the \textbf{mass} of $K$, which is a quantity strongly related to the class number but in many ways better behaved; this can be calculated via the special value of an appropriate $L$-function. 

\begin{definition}
The mass of a compact open subgroup $K$ of $G(\mathbb{A}_f)$, where $G$ is a connected reductive algebraic group compact at $\infty$, is the quantity
\[\Mass(K) = \sum_{[g] \in G(\mathbb{Q}) \backslash G(\mathbb{A}_f) / K} \frac{1}{\#\Gamma_g}.\]
\end{definition}

A mass formula for unitary groups is given in \cite{gan-hanke-yu}. This applies to a somewhat more general class of unitary groups than we consider here, so we shall state a special case: if $G$ is the unitary group $U_{n,E}$ as defined above and $n$ is odd, and $K = G(\hat{\mathbb{Z}})$, then
\[\Mass(K) = \frac{1}{2^n} \cdot L(M) \cdot \tau(G) \cdot \prod_{p} \lambda_p.\]
Here $L(M)$ is the value at $s = 0$ of the $L$-function of the motive of $G$; $\tau(G)$ is the Tamagawa number of $G$; and the local factors $\lambda_p$ are 1 unless $p$ is ramified in $E$, in which case they are $\frac{1}{2}$. (For $n$ even, there is a similar formula but the local factors are all 1.)

Formulae for $\tau(G)$ and $L(M)$ are known, and are given in \cite{gan-hanke-yu}. For unitary groups we have $\tau(G) = 2$; and if $\chi$ is the Dirichlet character corresponding to the quadratic extension $E / \mathbb{Q}$, then 
\[L(M) = \prod_{r = 1}^n L(1 - r, \chi^r).\]
which is easy to compute using generalised Bernoulli numbers as in chapter 4 of \cite{wash}.

\begin{remark}
Knowing the mass in advance allows one to determine whether the full class set has been found; but it is not completely clear that this will ever actually occur -- that is, that every $K$-class occurs in the decomposition of some Hecke operator. On the other hand, one can show (using the assumption that $E$ has unique factorisation and strong approximation for the derived subgroup $SU_{3,E}$, cf \cite[ch.~7]{platonov-rapinchuk}) that every $K$-class does have a representative supported at $p$ for every prime $p$. Hence if computing $T_{\mathfrak{p},1}$ and $T_{\mathfrak{p},2}$ for enough $\mathfrak{p}$ does not find enough $K$-classes, one can switch to computing multiple Hecke operators at some fixed small prime. However, I am not aware of any case in which this is necessary.
\end{remark}

\section{\texorpdfstring{Example: $E = \mathbb{Q}(\sqrt{-7})$}{Example: E = Q(sqrt(-7))}}

To illustrate the algorithm above, we shall carry out the calculations in the case $E = \mathbb{Q}(\sqrt{-7})$. This choice is convenient, as then $G$ is split at 2, and $\mathcal{O}_E^\times = \pm 1$, so $G(\mathbb{Q}) \cap K$ is small (it is the group of monomial matrices with entries in $\mathcal{O}_E^\times$, and thus has order 48). Let $\omega = \frac{1 + \sqrt{-7}}{2}$, so $\mathcal{O}_E = \mathbb{Z}[\omega]$.

We shall begin by calculating the mass of $K$. Here $\chi$ is the Kronecker symbol $\left(\frac{-7}{\bullet}\right)$; using standard methods we find that $L(0, \chi) = 1$, $L(-1, \chi^2) = -\frac{1}{12}$ and $L(-2, \chi^3) = -\frac{16}{7}$, so  $L(M) = \frac{4}{21}$. Hence
\[\Mass(K) = \frac{1}{8} \cdot \frac{4}{21} \cdot 2 \cdot \frac{1}{2} = \frac{1}{42}.\]

Since $\frac{1}{42} > \frac{1}{48}$, the mass of the principal $K$-class, there must be other $K$-classes; we shall find them by the bootstrapping argument explained above. 

If $p$ is any split prime and $\mathfrak{p}$ a choice of factor of $p$ in $E$, then we can apply Algorithm A to find the set $G(\mathbb{Q}) \cap K \eta_{\mathfrak{p},1} K$, taking $\lambda = \overline{\mathfrak{p}}$. For instance, let us take $p = 2$ and choose $\mathfrak{p}$ to be the ideal generated by $\omega$. Then we find 288 such matrices; since $G(\mathbb{Q}) \cap K$ has order 48, this gives us 6 right cosets $\gamma K$, which are represented by the three diagonal matrices with entries $1, 1, \frac{\omega}{1-\omega}$ in some order, and the three permutations of the rows of
\[\begin{pmatrix}
1 & 0 & 0 \\
0 & \frac{\omega}{2} & -\frac{\omega}{2} \\
0 & \frac{\omega}{2} & \frac{\omega}{2} 
\end{pmatrix}.\]

Lemma \ref{lemma:localbound} tells us that there is one coset missing, and the proof of the lemma constructs an element of this coset, namely the element $\alpha$ of $G(\mathbb{A}_f)$ that is 1 at all primes $p \ne 2$ and under the isomorphism $G(\mathbb{Q}_2) \cong \GL_3(\mathbb{Q}_2)$ corresponds to 
\[\begin{pmatrix}
1 & 0 & 0 \\
0 & 1 & 0 \\
-1 & -1 & 2 \\
\end{pmatrix}.\]

The associated group $\Gamma_\alpha$ has order $336$; indeed it is isomorphic to the direct product of $\pm 1$ with a simple subgroup isomorphic to $PSL_2(\mathbb{F}_7)$. Conveniently, this group contains $G(\mathbb{Z})$ as a subgroup of index 7; it is generated by $G(\mathbb{Z})$ and the element
\[
\tau = \begin{pmatrix}
-\frac{1}{2} & \frac{1 - \omega}{2} & \frac{1}{2} \\
-\frac{1}{2} & \frac{-1 + \omega}{2} & \frac{1}{2} \\
\frac{\omega}{2} & 0 & \frac{\omega}{2}\\
\end{pmatrix}.\]
of order 7.

Since $\frac{1}{48} + \frac{1}{336} = \frac{1}{42} = \Mass(K)$, that's all: the class number of $G$ is 2.

We can now complete the calculation of the Hecke operator $T_{\mathfrak{p},1}$ as follows. In the notation of 
the previous section, we are taking $\mu_1 = 1$ and $\mu_2 = \alpha$. We have decomposed $\mu_1 K \eta K$ as a union of 6 cosets of the form $\gamma \mu_1 K$ and 1 coset $\gamma \mu_2 K$. Searching for rational elements of $\mu_2 K \eta K$ (using Algorithm A with $\lambda = (1-\omega)^2$) finds 7 cosets already, so we do not need to consider the fourth possibility, which is whether there are any elements of $G(\mathbb{Q}) \cap \mu_2 K \eta K \mu_2^{-1}$. 

Since we have now found the complete class set, it is clear how to calculate $T_{\mathfrak{p}, i}$ for any split prime $p$. Note that we need to apply Algorithm A four times, with $\lambda = \overline{\mathfrak{p}}$, $\omega \overline{\mathfrak{p}}$, $\overline{\omega}\overline{\mathfrak{p}}$ and $2\overline{\mathfrak{p}}$, so we need to enumerate $r$-good matrices for $r = p$, $2p$ and $4p$.

\begin{remark} Calculating Hecke operators at the inert primes is also possible, but more awkward both theoretically and computationally. Since the maximal split torus in $G(\mathbb{Q}_p)$ for $p$ inert has rank 1, there is only one Hecke operator to deal with, but calculating its coset decomposition requires enumerating $r$-good matrices for $r = p^2$, which rapidly becomes computationally infeasible. Also, the eigenvalues of these Hecke operators are rather more indirectly related to Galois representations.
\end{remark}


\section{Experimental results: level 1}
\label{sect:weights}

In this chapter I'll summarise the calculations I have made for various small weights. Note that there are no forms of weight $(a,b,c)$ if $a + b + c$ is odd since the level group contains $-1$; and changing $c$ by a multiple of 2 just induces a twist, so I shall restrict to the case where $c=0$ if $a + b$ is even and $c = 1$ if $a + b$ is odd; this is what is meant by ``weight $(a,b)$''. 

\subsection{Dimensions of spaces}

These can be deduced from classical results of invariant theory for actions of groups on polynomial algebras. See table \ref{dimensions}.
\begin{table}[ht]
\centering\begin{tabular}{|c|c|cccccc|}
   \hline && & & a & & &\\
   \hline
   &   &        0 & 1 & 2 & 3 & 4 & 5 \\
   \hline
   & 0 &        2 & 0 & 1 & 0 & 3 & 0\\
   & 1 &        0 & 0 & 0 & 1 & 1 & 2\\
  b& 2 &        1 & 0 & 2 & 1 & 5 & 2\\
   & 3 &        0 & 1 & 1 & 4 & 3 & 6\\
   & 4 &        3 & 1 & 5 & 3 & 8 & 7\\
   & 5 &        0 & 2 & 2 & 6 & 7 & 10\\
   \hline
 \end{tabular}
\caption{\label{dimensions} Dimension of spaces of automorphic forms}
\end{table}

Note that if $a$, $b$ are small and $a + b$ is odd, then the spaces are often zero. Thus I will mostly restrict to the case of even weight.

\subsection{Forms arising from Gr\"ossencharacters}

For small weights, many of the forms that arise appear to have $T_{\mathfrak{p},1}$ eigenvalues of the form $\chi_1(\mathfrak{p}) + \chi_2(\mathfrak{p}) + \chi_3(\mathfrak{p})$, where $\chi_i(\mathfrak{p}) = \mathfrak{p}^r \overline{\mathfrak{p}}^s$ for some $(r,s)$. Note that there exists an unramified Gr\"ossencharacter of $E$ sending a uniformiser at $\mathfrak{p}$ to $\mathfrak{p}^r \overline{\mathfrak{p}}^s$ whenever $r + s$ is even, as $E$ has unique factorisation and unit group $\pm 1$; let us write $\chi(r,s)$ for this character.

In all cases other than the constant function, we find that the characters arising have $r + s = 2$. In table \ref{endoscopy1} below I list these, based on calculating Hecke operators at all split primes $p \le 60$. (I have assigned arbitrary numbers to $\Gal(\overline{E}/E)$-orbits of eigenforms, in order of discovery.)

\begin{table}[ht]
\centering
\begin{tabular}{|c|c|c|}
\hline
Weight  & Form  & Character\\
\hline
(0,0)   & 1     & $\chi(0,0) + \chi(1,1) + \chi(2,2)$ \\
(0,0)   & 2     & $\chi(2,0) + \chi(1,1) + \chi(0,2)$ \\
(2,0)   & 1     & $\chi(4, -2) + \chi(1,1) + \chi(0,2)$ \\
(4,0)   & 3     & $\chi(6,-4) + \chi(1,1) + \chi(0,2)$ \\
(2,2)   & 2     & $\chi(4,-2) + \chi(1,1) + \chi(-2,4)$ \\
\hline
\end{tabular}
\caption{\label{endoscopy1} Forms arising from Gr\"ossencharacters}
\end{table}

\subsection{Forms arising from classical modular forms}

Other forms of small weight appeared to have $T_{\mathfrak{p}, 1}$ eigenvalues given by $\chi(r_1,s_1) + \chi(r_2, s_2) a_p$ for $a_p$ the Hecke eigenvalue attached to a classical modular newform. Note that this class includes all of the forms above (some in multiple ways), since for each $r$ there exists a CM-type modular form of weight $r-1$ with $a_p = \chi(r,-r) + \chi(-r,r)$. The automorphic forms on $U(3)$ that are associated with non-CM-type forms are listed in Table \ref{endoscopy2} below.

\begin{table}[ht]
\centering
\begin{tabular}{|c|c|c|c|c|}
\hline
Weight  & Form  & Character & Modular weight & Modular label\\
\hline
(4,0)   & 1     & $\chi(1,1) + \chi(0,-4)  a_p$ & 7 & 7k7B[3]1 \\
(4,0)   & 2     & $\chi(0,2) + \chi(1,-4)  b_p$ & 6 & 49k6E1 \\
(3,1)   & 1     & $\chi(1,1) + \chi(-1,-3) a_p$ & 7 & 7k7B[3]1 \\
(2,2)   & 1     & $\chi(1,1) + \chi(-2,-2) a_p$ & 7 & 7k7B[3]1 \\
(3,3)   & 1     & $\chi(1,1) + \chi(-3,-3) c_p$ & 9 & 7k9B[3]1 \\
\hline
\end{tabular}
\caption{\label{endoscopy2} Forms arising from non-CM-type modular newforms}
\end{table}

In the third column of the table, I have reused letters when the same modular form occurs more than once. The labels in the final column are those used in William Stein's online tables.

Note that above $\chi(r,s)$ was defined only for $r + s$ even. If $r+s$ is odd, there is no unramified Gr\"ossencharacter sending a uniformiser at $\mathfrak{p}$ to $\mathfrak{p}^r \overline{\mathfrak{p}}^s$ -- indeed, this is not well-defined as the choice of generator of the ideal $\mathfrak{p}$ is only defined up to $\pm 1$ -- but there is a unique Gr\"ossencharacter which is unramified outside $\lambda = \sqrt{-7}$,  whose restriction to $\mathcal{O}(E_\lambda)$ is the order 2 character of $\mathbb{F}_7^\times$, and whose restriction to $E_{\infty} = \mathbb{C}$ is $z \mapsto z^{-r} \overline{z}^{-s}$; this character sends a uniformiser at $\mathfrak{p}$ to $ \left( \frac{\mathfrak{p}}{\lambda} \right) \mathfrak{p}^r \overline{\mathfrak{p}}^s$, which is independent of the choice of generator of $\mathfrak{p}$. This is what is meant by $\chi(1, -4)$ in the above table.

\subsection{Interpretation: Satake parameters}

Recall that the Satake parameters of an unramified representation $\pi_p$ of $\GL_n(\mathbb{Q}_p)$ are the eigenvalues of the semisimple conjugacy class that is the Langlands parameter of the representation $\pi_p$. These can be obtained as the roots of the Satake polynomial 
\[ \sum_{i=0}^n (-1)^i p^{i(i-1)/2} t_{\mathfrak{p}, i} X^{n-i} \]
where $t_{\mathfrak{p},i}$ is the eigenvalue of $T_{\mathfrak{p},i}$ on the $\GL_3(\mathbb{Z}_p)$-invariants of $\pi_p$; see \cite[\S I.8.4]{chenevier-thesis}. 

We introduce the following definition, from \cite[\S 6.4.5]{bellaiche-chenevier}: 

\begin{definition}
 Let $X$ be the multiset of Satake parameters (with multiplicities); then $\pi_p$ is is said to be {\it almost tempered} if there is a partition of $X$ into subsets $X_i$ such that each $X_i$ is of the form $\{y, y/q, \dots, y/q^{m_i-1}\}$, for $m_i = |X_i|$, and the real number $|\prod_{x \in X_i} x|^{1/m_i}$ does not depend on $i$; here $|\cdot|$ denotes the complex absolute value.
\end{definition}

If $\pi_p$ is tempered, which means in ths context that it is a subquotient of the parabolic induction of a unitary character of the standard Borel subgroup, then all of the Satake parameters have complex absolute value $p^{(n-1)/2}$, and so it is almost tempered with all the $X_i$ singletons; but not all almost tempered representations are of this type -- for example, a representation with Satake parameters $\{p^{1/2}, p^{-1/2}, 1\}$ would be almost tempered.\footnote{In general, the almost tempered representations are exactly those that are twists by unramified characters of parabolic inductions from a unitary character of the Levi subgroup of some parabolic; if this parabolic is not minimal, the representation will not be tempered. The above example corresponds to the trivial character of the Levi subgroup $\GL_2 \times \GL_1 \subset \GL_3$.} However, if $\pi_p$ is generic (admits a Whittaker model) and almost tempered, it must be tempered.

A conjecture due to Arthur predicts that if $\pi_p$ is the local factor of a cuspidal automorphic representation of any unitary group, then $\pi_p$ is almost tempered, which implies that if $\pi_p$ is generic, all the Satake parameters must have complex absolute value $p$. See \cite[\S 6.4]{bellaiche-chenevier} for further discussion.

Since our form of $U_{3,E}$ is compact at $\infty$, all automorphic forms are cuspidal. For all the forms we have encountered, the local factors at all split primes calculated are generic (except for the constant form, which is never generic). For the forms in Table \ref{endoscopy1} and Table \ref{endoscopy2} above, the Satake parameters are easy to read off: for Table \ref{endoscopy1}, the Satake parameters are the values at $p$ of the three characters listed, and for Table \ref{endoscopy2}, the Satake parameters for $\chi_1 + \chi_2 a_p$ are $\{ \chi_1(\mathfrak{p}), \chi_2(\mathfrak{p})\alpha, \chi_2(\mathfrak{p}) \beta\}$ where $\alpha$ and $\beta$ are the Satake parameters associated to the form $f$. The normalisation that is conventional for classical modular forms of weight $k$ gives both $\alpha$ and $\beta$ absolute value $p^{(k-1)/2}$. Thus the characters $\chi(r,s)$ in Table \ref{endoscopy1} are forced to have $r + s = 2$; and for characters of the form $\chi(r_1, s_1) + \chi(r_2,s_2) a_p$, we'd better have $r_1 + s_1 = 2$ and $r_2 + s_2 = 3-k$.

\subsection{Galois representations}

Let $f$ be an eigenform of weight $(a,b,c)$ and level $G(\hat{\mathbb{Z}})$. By results of Blasius and Rogawski \cite{blasius-rogawski}, for each prime $\ell$ there is a semi-simple representation
\[ \rho_f : \Gal(\overline{E}/E) \to \GL_3(\overline{\mathbb{Q}}_\ell)\]
unramified outside $7\ell$ and crystalline at $\ell$, such that for $p = \mathfrak{p} \overline{\mathfrak{p}}$ any other prime split in $E$ not dividing $7\ell$ the characteristic polynomial of $\rho_f(\Frob_{\mathfrak{p}})$ should be the Satake polynomial of $f$ at $p$ as defined above.

(The proofs of this statement in the literature assume that the weight is ``regular'', i.e.~$a \ne 0$ and $b \ne 0$. For non-regular weights, $\ell$-adic interpolation in the manner of \cite[Ch.~1]{chenevier-thesis} gives a representation with values in $\GL_3(\mathbb{C}_\ell)$ as long as $\ell$ splits in the field $E$, but it is not clear whether it is crystalline. It is, however, necessarily Hodge-Tate.)

Here $t_{\mathfrak{p},i}$ is the eigenvalue of $T_{\mathfrak{p},i}$ on $f$, so in particular $t_{\mathfrak{p},1}$ is the trace of Frobenius. Thus if the formulae conjectured above for the $T_{\mathfrak{p},1}$ eigenvalues of small weight forms are correct, the associated representations $\rho_f$ decompose as direct sums of smaller-dimensional representations (sums of three characters, or of a character and a twist of a 2-dimensional modular representation).

If $\ell$ is a prime split in $E$, then we can restrict the representation to a decomposition group at either of the primes above $\ell$. Since this local representation is crystalline, it is certainly Hodge-Tate; but the weights at these two different primes are not generally the same -- for one decomposition group the weights will be $-b+c, c+1, a + c + 2$ and for the other the weights will be $-a-c, -c + 1, b - c + 2$. This fits the tables above, as $\chi(a,b)$ has weight $a$ at one prime and $b$ at the other, and the representation associated to a modular form of weight $k$ has weights $0$ and $k-1$ (whichever prime is used). 

\subsection{Matching weights}

Given a fixed weight $(a,b,c) \in \mathbb{N} \times \mathbb{N} \times \mathbb{Z}$, the constraints of the valuations of Satake parameters and the Hodge-Tate weights at the two primes above $\ell$ leave few possibilities for endoscopic forms at this weight. 

Firstly, let us consider the case of automorphic forms $f$ arising from the endoscopic subgroup $U_1 \times U_1 \times U_1$. These correspond to Hecke eigenvalues that are the sum of three Gr\"ossencharacters of $E$, as in table \ref{endoscopy1}. If $f$ is to have level 1, then the three Gr\"ossencharacters must be unramified; so the $T_{\mathfrak{p}, 1}$ eigenvalues must be of the form $\chi(r_1,s_1) + \chi(r_2, s_2) + \chi(r_3, s_3)$ with $r_i+s_i$ even for each $i$.

We shall suppose $c = 0$. Then if the local factor of $f$ is tempered at some prime split in $E$, we must have $r_i + s_i = 2$ for each $i$. Then it is easy to see that the only possibility is $\chi(a+2, -a) + \chi(1,1) + \chi(-b, b + 2)$, which has level 1 if $a$ and $b$ are both even. This explains why we have seen one eigenform of this type when $a$, $b$ are both even and none otherwise.

Assuming Arthur's conjecture, if $f$ is not tempered, then it is almost tempered, two of its Satake parameters have ratio $p$, and we must have (up to a permutation) either $r_1 + s_1 = 4$, $r_2 + s_2 = 2$, and $r_3 + s_3 = 0$, or $r_1 + s_1 = 2$, $r_2 + s_2 = 3$, and $r_3 + s_3 = 1$. In the first case we also must have $r_2 = r_1 - 1$ and $r_3 = r_1 - 2$, and the same for the $s_i$; this can only occur for the constant form of weight $(0,0)$. The second case cannot occur at level 1, as we must have $r_i + s_i$ even.

The situation for endoscopic lifts from $U_2 \times U_1$ is more complicated. Let us assume our putative automorphic form $f$ on $U_3$ is generic at all places, has weight $(a,b,0)$ for integers $a,b \ge 0$, and has Hecke eigenvalues of the form $\chi(r_1, s_1) + \chi(r_2, s_2) a_p$ where $a_p$ is the Hecke eigenvalue of a classical modular cusp form $g$ of weight $k$. (Note that this includes the previous case, as endoscopic lifts from $U_1 \times U_1$ to $U_2$ correspond to modular forms $g$ of CM type.) The Galois representation associated to $g$ has Hodge-Tate weights $0$ and $k-1$, and its Satake parameters at $p$ have complex absolute values $p^{(k-1)/2}$, so we must have $r_1 + s_1 = 2$ and $r_2 + s_2 = 3-k$. Checking all the possible cases, we find that there are potentially 3 families of automorphic eigenforms on $U_3$, with $T_{\mathfrak{p},1}$ eigenvalues and weights given by
\begin{align*}
\chi(1, 1) + \chi(-s, 3-k+s)a_p &\text{ weight } (k-3-s, s, 0) \\
\chi(2-r, r) + \chi(1, 2-k)a_p &\text{ weight } (k-2,r-2,0) \\
\chi(r, 2-r) + \chi(2-k, 1)a_p &\text{ weight } (r-2,k-2,0). \\
\end{align*}

\subsection{A non-endoscopic form}

Weight $(3,3)$ is the first case where the action of the Hecke operators is not diagonalisable over $E$. One obtains two Hecke-invariant subspaces of dimension 2, corresponding to Galois orbits of eigenforms.

The first subspace splits over $E(\sqrt{46})$, and the $T_{\mathfrak{p},1}$-eigenvalues are $\chi(1,1) + \chi(-3,-3) c_p$, where $c_p$ is the unique Galois orbit of modular eigenforms of weight 9 for $\Gamma_1(7)$ such that $c_2$ is defined over this field. In this case, $T_{\mathfrak{p}, 2}$ coincides with $T_{\mathfrak{p},1}$.

The second space is more interesting. Here $T_{\mathfrak{p},1}$ and $T_{\mathfrak{p},2}$ commute with each other and have the same eigenvalues, but are not identical; so for the two eigenforms in the orbit, the $T_{\mathfrak{p},1}$-eigenvalue of one form is the $T_{\mathfrak{p},2}$-eigenvalue of the other form and vice versa. The eigenvalues are all defined over $\mathbb{Q}(\sqrt{-259})$, and the Satake polynomial for one of the primes above 2 is 
\[X^3 - \frac{-7 + \sqrt{-259}}{8} X^2 + \frac{-7 - \sqrt{-259}}{4} X - 8.\]

Since the Satake polynomial is irreducible over $E(\sqrt{-259})$, the Satake parameters are all $\Gal(\overline{E}/E)$-conjugate; so none of them arise from Gr\"ossencharacters of $E$, and thus the corresponding $\ell$-adic Galois representations must be irreducible of dimension 3.

Let's check the valuations of the Satake parameters. It can be shown \cite{DL-cubics} that a cubic over $\mathbb{C}$ has all its roots of equal valuation if and only if, after scaling so the constant term is 1, it is of the form $T^3 + \overline{a} T^2 + a T + 1$ with $a$ lying inside a deltoid in the complex plane with centre 0 and one cusp at 3; this is equivalent to $a = x + iy$ with $x,y$ real such that
\[C(x,y) = (x^2 + y^2)^2 - 8x(x^2 - 3y^2) + 18(x^2 + y^2) - 27 \le 0.\]
In this case we find that $C(x,y) = -\frac{56727}{4096}$, so the condition is satisfied.


\section{Some higher level examples}
\label{sect:higher}

We shall now introduce some level structure at the prime 2. In preparation for future work in which I intend to consider overconvergent automorphic forms, we shall use a slightly different definition of automorphic forms in which the action is twisted to the other side. We shall fix a prime $p$ which is split in $E$ and an irreducible algebraic representation $V$ of $G$ (which will be defined over $\mathbb{Q}_p$), and consider functions 
\[f : G(\mathbb{A}_f) \to V\quad\text{such that}\quad f(\gamma g k) = k_p^{-1} f(g)\]
for all $\gamma \in G(\mathbb{Q})$ and $k \in K$ for some open compact $K$. This is isomorphic to the space considered above, via the map $f(g) \mapsto g_p \circ f(g)$ (if we fix a choice of embedding of $E \hookrightarrow \mathbb{Q}_p$), but has the advantage that we only ever need to consider the action on $V$ of elements of $K_p$; this is useful when considering $p$-adic interpolation, as in \cite{chenevier-thesis} or my University of London PhD thesis (in preparation).

The action of the Hecke algebra on this space is slightly different: we define an action of a double coset $K g K$ on automorphic forms of level $K$ by writing $K g K = \bigsqcup g_{i} K$, and defining
\[ ([K gK] \circ f)(x) = \sum_{i} (g_{i})_p f(x g_i).\]
It's easily verified that this gives a well-defined linear endomorphism of the space of automorphic forms. 

\subsection{\texorpdfstring{$L$-classes}{L-classes}}

We shall set $p=2$; and we will define our level group $L$ to be $G(\mathbb{Z}_\ell)$ at all $\ell \ne 2$, and at 2 the subgroup $L_p$ of $\GL_3(\mathbb{Z}_2)$ which reduces mod 2 to
\[\begin{pmatrix} * & * & * \\ 0 & * & * \\ 0 & * & *\end{pmatrix}.\]
This is the parahoric subgroup associated to a parabolic subgroup of $\GL_3$. We now need to calculate $L$-classes and decompositions of Hecke operators for this case. We shall do this by comparing $L$ to the group $K = \GL_3(\hat{\mathbb{Z}})$ considered above, using the following lemma:

\begin{lemma}
Let $U' \subset U$ be any finite index subgroups of $G(\mathbb{A}_f)$. Then we have $\Mass(U') = [U : U'] \Mass(U)$.
\end{lemma}

\begin{proof} This is an essentially trivial but fiddly manipulation.

\begin{align*}
\Mass(U') &= \sum_{[v] \in G(\mathbb{Q}) \backslash G(\mathbb{A}_f) / U'} \frac{1}{|G(\mathbb{Q}) \cap v U' v^{-1}|}\\
&= \sum_{[\mu] \in G(\mathbb{Q}) \backslash G(\mathbb{A}_f) / U} \sum_{k \in U/U'} \frac{1}{|G(\mathbb{Q}) \cap \mu k U' k^{-1}\mu^{-1}|} \cdot \frac{1}{\text{size of orbit of $k$}}
\end{align*}
(where we consider the action of $U \cap \mu^{-1}G(\mathbb{Q})\mu$ on $U / U'$ by left multiplication)
\begin{align*}
&= \sum_{[\mu] \in G(\mathbb{Q}) \backslash G(\mathbb{A}_f) / U} \sum_{k \in U/U'} \frac{1}{|G(\mathbb{Q}) \cap \mu k U' k^{-1}\mu^{-1}|} \cdot \frac{|G(\mathbb{Q}) \cap \mu k U' k^{-1}\mu^{-1}|}{|G(\mathbb{Q}) \cap \mu U \mu^{-1}|}\\
&= \Mass(U) \cdot \left(\sum_{k \in U/U'} 1\right)\\
&= \Mass(U) \cdot [U : U']
\end{align*}
\end{proof}

It's readily seen that $L_p$ is the subgroup of matrices in $\GL_3(\mathbb{Z}_2)$ whose reduction mod 2 stabilises $(*, 0, 0)^T$ under the left multiplication action on column vectors; $\GL_3(\mathbb{Z}_2)$ clearly acts transitively on these, so $[K : L] = |\mathbb{P}^2(\mathbb{F}_2)| = 7$. Hence $\Mass(L) = \frac{7}{42} = \frac{1}{6}$.

The proof of the above lemma indicates how to work out the $L$-classes, by calculating, for each $K$-class $G(\mathbb{Q}) \mu K$, how $K \cap \mu^{-1}G(\mathbb{Q})\mu$ acts by left multiplication on $K/L$, which we can identify with $\mathbb{P}^2(\mathbb{F}_2)$. 

For the identity $K$-class, the group $G(\mathbb{Q}) \cap K$ consists of a product of diagonal $\pm 1$'s (which all reduce to the identity and can be ignored) and permutation matrices. So the action corresponds to permutation of coordinates on $\mathbb{P}^2(\mathbb{F}_2)$, which has exactly 3 orbits classified by the number of nonzero terms. Coset representatives are given by: the identity (having stabiliser of order 16); the matrix $\begin{pmatrix} 1 & 0 & 0 \\ 1 & 1 & 0 \\ 0 & 0 & 1\end{pmatrix}$ (again with 16); and the matrix $\begin{pmatrix} 1 & 0 & 0 \\ 1 & 1 & 0 \\ 1 & 0 & 1 \end{pmatrix}$ with stabiliser of order 48. 

For the nonidentity $K$-class, the extra generator $\tau$ of $G(\mathbb{Q}) \cap \alpha K \alpha^{-1}$ conjugates to an element of $G(\mathbb{A}_f)$ whose matrix reduces mod $p$ to 
\[\begin{pmatrix} 1 & 1 & 1 \\ 1 & 0 & 1 \\ 1 & 0 & 0\end{pmatrix}\]
This acts transitively on the elements of $\mathbb{P}^2(\mathbb{F}_2)$, so $G(\mathbb{Q}) \alpha K = G(\mathbb{Q}) \alpha L$. The stabilizer group must thus have size 48, and it is easy to enumerate its elements as it is a subgroup of the finite group $\Gamma_{\alpha}$ of the previous two chapters.

\subsection{\texorpdfstring{The operator $U$}{The operator U}}

Motivated again by the theory of overconvergent automorphic forms, we consider the action of the Hecke operator corresponding to the double coset $L \eta_{\mathfrak{p}, 2} L$, where $\mathfrak{p}$ is the prime dividing 2 generated by $\omega$. 
As the role it plays is roughly analogous to the Atkin-Lehner $U_p$-operator in Coleman's theory, we shall denote it by the letter $U$.

To decompose $L \eta_{\mathfrak{p}, 2} L$ into single cosets locally, we observe that if $\eta = \eta_{\mathfrak{p}, 2}$, then $L_p \cap \eta L_p \eta^{-1}$ consists of those elements of $L_p$ whose first column reduces to $(*, 0, 0)^{T} \bmod 4$. This is an index 4 subgroup of $L_p$ and coset representatives $\phi_i$ may be obtained by taking the first column to be $(1,0,0)^T$, $(1,2,0)^T$, $(1,0,2)^T$ and $(1,2,2)^T$. We hence obtain a decomposition of $L\eta L$ as $\bigsqcup_{i=1}^4 \phi_i \eta L$.

As in the level $K$ case, the final step is to write each single coset in terms of our chosen $L$-class representatives $\mu_1 \dots \mu_4$, and more generally to do the same for the products $\mu_r \phi_s \eta L$. This is easy using an adaptation of Algorithm A: we note that if $\gamma \in G(\mathbb{Q})$ satisfies $\gamma \mu_t \in \mu_r \phi_s \eta L$, then the largest denominator that can possibly occur in $\gamma$ has norm either $2$, $4$ or $8$, depending on $(r,s,t)$. We then generate lists of the $2$, $4$ and $8$-good matrices, and for each coset $\mu_r \phi_s \eta L$, we use these to find a representative for this coset in the form $\gamma \mu_t$ (testing each of the possible values of $t$ in turn until one works). Note that we must explicitly test $\phi_s^{-1}\mu_r^{-1} \gamma \mu_t^{-1}$ for membership of $L$, rather than dodging this check using elementary factors as before. The results of this computation are shown in Table \ref{uoperator}.

\begin{remark}
Since the group $L_p$ we have chosen is a parahoric subgroup, one can attach invariants to double cosets $L_p g L_p \subset \GL_3(\mathbb{Q}_p)$ (analogous to elementary factors) using the Bruhat--Iwahori decomposition. One can probably use this to calculate the Hecke operators more efficiently; but the slower method above is more flexible and would apply to any group $L$ for which we could compute $K / L$ and $L / (L \cap \eta L \eta^{-1})$.
\end{remark}

\begin{table}[ht]
\begin{align*}
\mu_1 \phi_1\eta L &= \begin{pmatrix} 1 & 0 & 0 \\ 0 & 1 - \frac{w}{2}& 0 \\0 & 0 & 1 - \frac{w}{2} \end{pmatrix} \mu_1 L &
\mu_3 \phi_1 \eta L &= \begin{pmatrix}\frac{w}{2} & -\frac{w}{2} & 0 \\ \frac{w}{2} & \frac{w}{2} & 0 \\ 0 & 0 & 1 \end{pmatrix} \mu_4 L \\
\mu_1 \phi_2\eta L &= \begin{pmatrix} 1 & 0 & 0 \\ 0 & 1 - \frac{w}{2}& 0 \\0 & 0 & 1 - \frac{w}{2} \end{pmatrix} \mu_2 L &
\mu_3 \phi_2 \eta L &= \begin{pmatrix}\frac{w}{2} & \frac{w}{2} & 0 \\ -\frac{w}{2} & \frac{w}{2} & 0 \\ 0 & 0 & 1 \end{pmatrix} \mu_4 L\\
\mu_1 \phi_3\eta L &= \begin{pmatrix} 1 & 0 & 0 \\ 0 & 0 & 1 - \frac{w}{2} \\0 & 1 - \frac{w}{2} & 0 \end{pmatrix} \mu_2 L &
\mu_3 \phi_3 \eta L &= \begin{pmatrix} 1 & 0 & 0 \\0 & \frac{w}{2} & \frac{w}{2}\\0 & \frac{w}{2} & -\frac{w}{2} \end{pmatrix} \mu_4 L \\
\mu_1 \phi_4 \eta L &= \begin{pmatrix}1 & 0 & 0 \\ 0 & 1 - \frac{w}{2}& 0 \\0 & 0 & 1 - \frac{w}{2} \end{pmatrix} \mu_3 L &
\mu_3 \phi_4 \eta L &= \begin{pmatrix} -\frac{w}{2} & \frac{w}{2} & 0 \\ \frac{w}{2} & \frac{w}{2} & 0 \\ 0 & 0 & 1\end{pmatrix} \mu_4 L \\
\mu_2 \phi_1 \eta L &= \begin{pmatrix}\frac{w}{2} & -\frac{w}{2} & 0 \\ \frac{w}{2} & \frac{w}{2} & 0 \\ 0 & 0 & 1-\frac{w}{2} \end{pmatrix} \mu_1 L &
\mu_4 \phi_1 \eta L &= \begin{pmatrix}-\frac{w}{2} & -\frac{w}{4} + \frac{1}{2} & -\frac{w}{4} + \frac{1}{2} \\ 0 & -\frac{w}{4} - \frac{1}{2} & \frac{w}{4} + \frac{1}{2} \\ \frac{w}{2} & -\frac{w}{4} + \frac{1}{2} & -\frac{w}{4} + \frac{1}{2}\end{pmatrix} \mu_1 L \\
\mu_2 \phi_2 \eta L &= \begin{pmatrix}\frac{w}{2} & \frac{w}{2} & 0 \\ -\frac{w}{2} & \frac{w}{2} & 0 \\ 0 & 0 & 1-\frac{w}{2} \end{pmatrix} \mu_1 L &
\mu_4 \phi_2 \eta L &= \begin{pmatrix}-\frac{w}{4} + \frac{1}{2} & -\frac{w}{2} & -\frac{w}{4} + \frac{1}{2} \\ \frac{w}{4} + \frac{1}{2} & 0 & -\frac{w}{4} - \frac{1}{2} \\ -\frac{w}{4} + \frac{1}{2} & \frac{w}{2} & -\frac{w}{4} + \frac{1}{2}\end{pmatrix} \mu_1 L \\
\mu_2 \phi_3 \eta L &= \begin{pmatrix}0 & \frac{w}{2} & -\frac{w}{2} \\ 0 & \frac{w}{2} & \frac{w}{2} \\ 1-\frac{w}{2} & 0 & 0 \end{pmatrix} \mu_2 L & 
\mu_4 \phi_3 \eta L &= \begin{pmatrix}-\frac{w}{2} & -\frac{w}{4} + \frac{1}{2} & -\frac{w}{4} + \frac{1}{2} \\ 0 & -\frac{w}{4} - \frac{1}{2} & \frac{w}{4} + \frac{1}{2} \\ \frac{w}{2} & -\frac{w}{4} + \frac{1}{2} & -\frac{w}{4} + \frac{1}{2} \end{pmatrix}  \mu_3 L\\
\mu_2 \phi_4 \eta L &= \begin{pmatrix}0 & \frac{w}{2} & \frac{w}{2} \\ 0 & -\frac{w}{2} & \frac{w}{2} \\ 1-\frac{w}{2} & 0 & 0 \end{pmatrix} \mu_2 L &
\mu_4 \phi_4 \eta L &= \begin{pmatrix} -\frac{w}{4} + \frac{1}{2} & -\frac{w}{2} & -\frac{w}{4} + \frac{1}{2} \\ -\frac{w}{4} - \frac{1}{2} & 0 & \frac{w}{4} + \frac{1}{2} \\ -\frac{w}{4} + \frac{1}{2} & \frac{w}{2} & -\frac{w}{4} + \frac{1}{2} \end{pmatrix} \mu_1 L
\end{align*}
\caption{\label{uoperator} Decomposition of the double coset $L \eta L$}
\end{table}

\subsection{Automorphic forms and slopes}

The coset representatives calculated in Table \ref{uoperator} above can now be used to calculate automorphic forms of level $L$ using the same machinery as before. The choice of level structure and of $U$ operator should imply that the resulting objects interpolate well 2-adically as the $a$ component of the weight varies. See Table \ref{slopes} for a table of the 10 smallest 2-adic slopes of the $U$ operator on forms of weight $(a,0,0)$, for various integers $a$ (necessarily even). This table shows clear evidence of the $2$-adic local constancy of the slopes.

\begin{table}
\centering
\begin{tabular}{|c|c c c c c c c c c c|}
\hline
$a$ & \multicolumn{10}{c|}{Slopes}\\
\hline
0 &     0 & 1 & 2 & 2 & & & & & &\\
2 &     0 & 1 & $\frac{3}{2}$ & $\frac{3}{2}$  & 3 & 4 & & & &\\
4 &     0 & 1 & 2 & $\frac{5}{2}$ & $\frac{5}{2}$ & 3 & 3 & 4 & 5 & 5\\
6 &     0 & 1 & $\frac{3}{2}$ & $\frac{3}{2}$ & 3& 3& $\frac{7}{2}$ & $\frac{7}{2}$ & $\frac{11}{2}$& $\frac{11}{2}$\\
8 &     0 & 1 & 2 & 2 & 2 & 4 & 4 & $\frac{9}{2}$ & $\frac{9}{2}$& $\frac{9}{2}$\\
10 &    0 & 1 & $\frac{3}{2}$ & $\frac{3}{2}$&  3& 5 &5 &5 &5 & $\frac{11}{2}$ \\
12 &    0 & 1 &2& $\frac{8}{3}$& $\frac{8}{3}$& $\frac{8}{3}$& 3& 4& 5& 6\\
14 &    0 &1& $\frac{3}{2}$& $\frac{3}{2}$& 3 &3 &$\frac{7}{2}$ &$\frac{7}{2}$ &$\frac{11}{2}$ &$\frac{11}{2}$ \\
16 &    0 &1 &2 &2 &2 &4 &4 &$\frac{14}{3}$ &$\frac{14}{3}$ & $\frac{14}{3}$\\
\hline
\end{tabular}
\caption{\label{slopes} The 10 smallest slopes of $U$ acting on forms of level $L$ and weight $(a,0,0)$, for even $a \le 16$.}
\end{table}

Note that the slopes of $U$ are not remotely locally constant with regard to variation in the $b$ direction in the weight lattice; indeed the dimension of its ordinary subspace grows as $b$ gets large (in the archimedean sense!). This is unsurprising, as $U$ does not define a compact endomorphism on the space of overconvergent forms as defined in \cite{chenevier-thesis}, but it is compact restricted to the smaller space of ``semi-classical'' automorphic forms associated to the parabolic subgroup corresponding to the parahoric $L_p$ -- see my forthcoming University of London PhD thesis for the definitions and properties of these spaces. This space corresponds to allowing $p$-adic variation in the $a$ direction but not the $b$ direction in the weight lattice. 

A striking aspect of these tables of slopes is the fact that they are often non-integral. This is very much contrary to the case of $\GL_2$ where slopes are ``usually'' integral, at least for weights in the centre of weight space (Wan has conjectured that for any fixed weight $k$ the denominators of the slopes of overconvergent forms at weight $k$ is bounded). Extending the computations above gave an example of a form with slope $\frac{74}{7}$ at weight $(5,5,0)$.

\section{Implementation details}

\subsection{Algorithms for \texorpdfstring{$r$}{r}-good matrices}
\label{sect:rgood}

The time-consuming step in Algorithm A, which was used in several places, is to calculate all elements of $M_3(\mathcal{O}_E)$ satisfying $m m^\dagger = r$; to calculate the Hecke operator at a split prime $p$ one needs to consider $r = p$, $2p$, and $4p$, so to have a reasonable supply of primes one needs to consider $r \approx 100$. Clearly this is equivalent to finding all orthogonal triples of vectors of length $r$ in $\mathcal{O}_E^3$, regarded as a 6-dimensional $\mathbb{Z}$-lattice. I experimented with two different algorithms, both of which have the same first step of generating a complete list $\mathcal{L}$ of the vectors of length $r$. This list is typically quite large -- since the theta-function of the lattice $\mathcal{O}_E^3$ is a modular form of weight 3 and level $\Gamma_1(7)$, it is easily seen that for $r = p, 2p, 4p$ with $p$ a prime split in $E$, there are $\Theta(p^2)$ such vectors (with leading terms that can be calculated precisely using Eisenstein series in each case).

The obvious algorithm is to regard this as a graph triangle enumeration problem: consider the graph whose vertices are the vectors in $\mathcal{L}$ and where vertices $u,v$ are joined by an edge if they are orthogonal. This graph is typically not very dense, so from an appropriate sparse representation of its adjacency matrix one can enumerate all triangles in it quickly by a vertex-iterator approach. The most computationally difficult step is calculating the adjacency matrix: testing all possible pairs of elements of $\mathcal{L}$ will take $O(p^4)$ steps.

A more devious alternative is to calculate, for each vector $v$ of length $r$, a basis for the orthogonal complement of $v$ in $\mathcal{O}_E^3$, as a 4-dimensional $\mathbb{Z}$-lattice; enumerate vectors of length $r$ in this; and find the orthogonal complement of these, to give orthogonal triples. This option turns out to be very much quicker for large $r$, as the orthogonal complement lattices typically have very few short vectors.

\subsection{Programs included with this paper}

A selection of computer programs that were used in the calculations above can be downloaded from my personal website at \url{http://www.ma.ic.ac.uk/~dl505/maths/programs/programs.html}. These use a combination of the Sage \cite{sage} and Magma \cite{magma} computer algebra systems; users who do not have access to Magma can still use the code for level 1 automorphic forms, but will be restricted to those primes for which precalculated lists of $r$-good matrices are included.

\section{Acknowledgements}

I would like to thank Gaetan Chenevier, for initially raising the question of computing automorphic forms on unitary groups and for many helpful discussions since, and my PhD supervisor Kevin Buzzard, for his unfailing encouragement and support. 

I am also grateful to Michael Stoll, for suggesting to me the second, more efficient lattice-based algorithm for $r$-good matrices described in \S \ref{sect:rgood} above.

\providecommand{\bysame}{\leavevmode\hbox to3em{\hrulefill}\thinspace}
\renewcommand{\MR}[1]{\MRhref{#1}{MR #1.}}
\renewcommand{\MRhref}[2]{%
  \href{http://www.ams.org/mathscinet-getitem?mr=#1}{#2}
}

\end{document}